\documentclass[11pt] {article} 
\usepackage[backref]{hyperref}
  \usepackage{amsmath}
    \usepackage{amssymb}
    
    \usepackage{pdfsync}

\newtheorem{example}{Example}[section]
\newtheorem{theorem}{Theorem}[section]
\newtheorem{lemma}{Lemma}[section]
\newtheorem{corollary}{Corollary}[section]

\newtheorem{remark}{Remark}[section]

\newcommand{\eqnsection}{
   \renewcommand{\theequation}{\thesection.\arabic{equation}}
   \makeatletter
   \csname @addtoreset\endcsname{equation}{section} 
   \makeatother}

\def \ov{\overline}

\def \be{\begin{equation}}
\def \ee{\end{equation}}
\def \bt{\begin{theorem}} 
\def \et{\end{theorem}}
\def \bl{\begin{lemma}} 
\def \el{\end{lemma}}
\def \bea{\begin{eqnarray}}
\def \eea{\end{eqnarray}}
\def \bas{\begin{eqnarray*}}
\def \eas{\end{eqnarray*}}

\def \al{\alpha}
\def \bb{\beta}
\def \ga{\gamma}
\def \Ga{\Gamma}
\def \de{\delta}

\def \ep{\epsilon}
\def \vep{\varepsilon}
\def \la{\lambda}

\def \si{\sigma}

\def \ze{\zeta}

\def \ff{\infty}
\def \wh{\widehat}
\def \wt{\widetilde}

\def \cd{\,\cdot\,}

\def\stl{\stackrel{law}{=}}

\def \AA{{\cal A}}
\def \BB{{\cal B}}

\def \HH{{\cal H}}
\def \II{{\cal I}}

\def \RR{{\cal R}}

\def \({\left(}
\def \){\right)}

\def \nn{\nonumber}
\def \Proof{\noindent{\bf Proof $\,$ }}

\def \bc{\begin{center} }
\def \ec{\end{center} }
\def \bs{\begin{slide} }
\def \es{\end{slide} }

\def\square{{\vcenter{\vbox{\hrule height.3pt
        \hbox{\vrule width.3pt height5pt \kern5pt
           \vrule width.3pt}
        \hrule height.3pt}}}}
\def\qed{{\hfill $\square$ \bigskip}}

\eqnsection

\begin{document}

\title{  Conditions for permanental processes to be unbounded}

  \author{  Michael B. Marcus\,\, \,\, Jay Rosen \thanks{Research of   M. B. Marcus and J. Rosen 
   was partially supported by  grants from the National Science Foundation and
 PSC CUNY.   }}
\maketitle

\begin{abstract}   An $\al$-permanental process $\{X_{ t},t\in T \}$ is a stochastic process determined by   a kernel $K=\{K(s,t),s,t\in T \}$, with the property that for all $t_{1},\ldots,t_{n}\in T $,   $ |I+K( t_{1},\ldots,t_{n}) S|^{- \al}  $ is the Laplace transform of $(X_{t_{1}},\ldots,X_{t_{n}})$, where $ K( t_{1},\ldots,t_{n})$   denotes the matrix $\{K(t_{i}, t_{j})\}_{i,j=1}^{n}$ and $S$ is the diagonal matrix with entries $s_{1},\ldots,s_{n}  $. $ (X_{t_{1}},\ldots,X_{t_{n}})$ is called a permanental vector. 

Under the condition that $K$ is the potential density of a transient Markov process, 
 $(X_{t_{1}},\ldots,X_{t_{n}})$ is represented as a   random mixture of 
 $n$-dimensional random variables with components that are independent gamma  random variables.   This representation leads to  a Sudakov type inequality for the sup-norm of   $(X_{t_{1}},\ldots,X_{t_{n}})$ that is used to obtain   sufficient conditions for  a large class of permanental processes to be unbounded almost surely.   These results are used to obtain conditions for   permanental  processes associated with certain L\'evy processes to be unbounded. 
 
   Because $K$ is the potential density of a transient Markov process, for all $t_{1},\ldots,t_{n}\in T $, $A( t_{1},\ldots,t_{n}):= (K( t_{1},\ldots,t_{n}))^{-1}$ are $M$-matrices.     The results in this paper are obtained by working with these $M$-matrices.

  \end{abstract}

\bibliographystyle{amsplain}

\section{Introduction}\label{sec-1} An $R^{n}$ valued   $\al$-permanental random variable $X=(X_{1},\ldots, X_{n})$ is  a random variable with Laplace transform 
\begin{equation}
   E\(e^{-\sum_{i=1}^{n}s_{i}X_{i}}\) 
 = \frac{1}{ |I+KS|^{ \al}}   \label{int.1} 
 \end{equation}
for some $n\times n$ matrix $K$   and diagonal matrix $S$ with entries $s_{i}$, $1\le i\le n$, and $\al>0$. Permanental random variables were introduced by Vere-Jones, \cite{VJ}, who referred to them as random variables with  multivariate gamma distributions. (Actually he considered the moment generating function.) 

   \label{pp}   An $\al$-permanental process $\{X_{t},t\in T \}$ is a stochastic process  which has finite dimensional  distributions  that are $\al$-permanental vectors. The permanental process is determined by   a kernel $\{K(s,t),s,t\in T \}$, with the property that for all $t_{1},\ldots,t_{n}$ in $T $,  $\{K(t_{i},t_{j}),i,j\in [0,n] \}$ determines an $\al$-permanental random variable by (\ref{int.1}).  (Sometimes we  refer to these processes  simply as permanental processes.) Vere-Jones    briefly considers permanental processes in \cite{VJ}.  Note that when (\ref{int.1})  holds  for a kernel  $K(s,t)$ for all $\al>0$, the family of permanental processes obtained are infinitely divisible.  The permanental processes considered in this paper have this property.  
   
 Local times of Markov processes with symmetric potential densities   are related by isomorphism theorems to the squares of Gaussian processes.  Note that when $K$ is symmetric and positive definite and $\al=1/2$,  $ (\eta_{1}^{2}/2,\ldots,\eta_{n}^{2}/2)$,     where  $(\eta_{1}, \ldots,\eta_{n})$ is an $n$-dimensional   normal random variable  with mean zero and  covariance matrix $K$, is a 1/2-permanental process.   When $\al\ne 1/2$ or  $K$ is not symmetric, the isomorphism theorems can be generalized,   by replacing  the squares of the Gaussian processes by other permanental processes,   so that they    also hold for Markov processes with potential densities that are not symmetric.   To apply these  isomorphism theorems it is important to know sample path properties of permanental processes. 
   
 	In this paper we give a concrete representation of permanental vectors that is used to  obtain a Sudakov type inequality that gives  lower bounds for permanental  processes that only requires that the inverses of the matrices  $\{K(t_{i},t_{j}),i,j\in [0,n] \}$ are $M$-matrices. (It does not require that the matrices   are symmetric.)

  Since the definition of permanental processes requires that their finite dimensional distributions are permanental random  variables, 
a fundamental question is: For which matrices $K$   do there exist random variables $X$ satisfying (\ref{int.1})?  Vere-Jones answers this question but with criteria that are, in general,   very difficult  to verify.   On the other hand,   as we just pointed out, when $K$ is symmetric and positive definite and $\al=1/2$ then $X=(\eta_{1}^{2}/2,\ldots,\eta_{n}^{2}/2)$,     where  $(\eta_{1}, \ldots,\eta_{n})$ is an $n$-dimensional   normal random variable  with mean zero and  covariance matrix $K$.  

There are other cases in which   it is easy to see that the right-hand side of (\ref{int.1}) is the Laplace transform of an  $R^{n}$ valued  random variable. Recall that a gamma random variable is one with
 probability density function 
 \be
f(u,v;x) = \frac{v^{u} x^{u-1} e^{-v x}}{\Gamma(u)} \quad \text{ for } x \geq 0 \text{ and } u, v > 0,\label{1.2q}
 \ee and equal to 0 for $x\le0$, where
$\Gamma(u) =\int_{0}^{\ff}  x^{u-1} e^{-  x}\,dx 
$ is the gamma function.      The parameter $u$ is called  the shape of the gamma distribution and the parameter $v$ is called  the scale of the gamma distribution.

\medskip	  In this paper we describe a large class of infinitely divisible permanental random variables. We use $\xi_{u,v}$ to denote a random variable with  probability density function $f(u,v;x)$. 
The Laplace transform of $\xi_{u,v}$ is  
 \be 
    \int_{0}^{\ff}\frac{v^{u} x^{u-1} e^{-(v +s)x}}{\Gamma(u)}\,dx\\
   =\frac{1}{\(1+ \displaystyle \frac{s}{v} \)^{u}} =\frac{v^{u}}{\(v+ s \)^{u}}.\label{1.3q}
   \ee
  Therefore if $K$ is a diagonal matrix with entries $1/v_{i}$, (\ref{int.1})  is the Laplace transform of  $(\xi_{\al, v_{1}},\ldots,\xi_{\al, v_{n}})$, in which all the components are independent. Consequently, when  the right-hand side of (\ref{int.1})   is the Laplace transform of an $R^{n}$ valued  random variable  $X$, it is reasonable to say that $X$ has  an $n$-dimensional gamma distribution.  
   
 \medskip	   We assume that $|K|>0$. Therefore, $A=K^{-1}$ exists and we can also define $X$ by
\begin{equation}
   E\(e^{-\sum_{i=1}^{n}s_{i}X_{i}}\) 
 =     \frac{|A|^{\al}}{|A+S|^{\al}} \label{2.1}. 
 \end{equation}
   It turns out that it is simpler to  describe the random variables $X$ that are defined by matrices $K$   as in (\ref{int.1}), by focusing on $A$, and describing  the random variables $X$ that are defined by matrices $A$   as in (\ref{2.1}).

\medskip	 The results in this paper all depend on  a concrete representation of permanental random variables which we can obtain when the matrix $A$ in (\ref{2.1}) is a non-singular $M$-matrix.

\medskip	    Let  $C=\{ c_{ i,j}\}_{ 1\leq i,j\leq n}$ be an
$n\times n$ matrix. We call $C$ a positive matrix
and write
$C\geq 0$ if
$c_{ i,j}\geq 0$ for all
$i,j$.  

\medskip The matrix
$A$  is said to be a nonsingular
$M$-matrix   if
\begin{enumerate}
\item[(1)] $a_{ i,j}\leq 0$ for all $i\neq j$.
\item[(2)] $A$ is nonsingular and $A^{ -1}\geq 0$.
\end{enumerate}

    Theorem \ref{theo-9.1}  gives a representation of   $\al$ permanental vectors. It is rather technical and requires some preparation so we hold off presenting it until Section \ref{proofs}.    The following  consequence of Theorem \ref{theo-9.1} is our key to obtaining conditions for the paths of permanental  processes to be unbounded.

  \begin{theorem}\label{theo-cup} Let $X=(X_{1}, \dots, X_{n})$ be an $\al$-permanental vector with non-singular kernel $K$. Assume that $A=K^{-1}$ is an  $M$-matrix with diagonal entries $(a_{1}, \ldots, a_{n})$.  Then there exists a coupling between $X$ and an $n$-tuple $ \(\xi^{(1)}_{ \al,1},\ldots,\xi^{(n)}_{ \al,1}\)$ of independent identically distributed  copies of $\xi_{ \al,1}$ such that
     \begin{equation}
 X\geq   \(a_{1}^{-1}\xi^{(1)}_{ \al,1},\ldots,a_{n}^{-1}\xi^{(n)}_{ \al,1}\),\qquad a.s.\label{1.5qa}
  \end{equation} 
      \end{theorem}
 
This immediately implies the next theorem.

\bt\label{theo-9.2}	Let $X$ be as in Theorem \ref{theo-cup}.
Then    if  $f$ is an increasing function on  $R^{n}_{+}$   
 \begin{equation}
  E(f(X))\ge   E ( f(a^{-1}_{1}\xi^{(1)}_{ \al,1},\ldots,a^{-1}_{n}\xi^{(n)}_{ \al,1} ) ).\label{10.14}
  \end{equation}
Equivalently,
   \begin{equation}
  E f((a_{1}X_{1},\ldots,a_{n}X_{n}))\ge  E ( f( \xi^{(1)}_{ \al,1},\ldots, \xi^{(n)}_{ \al,1} ) ).\label{10.14z}
  \end{equation}
  \et

We call    (\ref{10.14}) the Permanental Inequality. We explain  in   Section  \ref{sec-10} that it is a generalization, in  a certain sense, of the Sudakov Inequality.

 \medskip	 
  It is shown in \cite{EK} that when $\{u(s,t),s,t\in T\}$ is the potential density of a transient Markov process with state space $T$, then   for any $\al>0$ , there exists an $\al$-permanental process with kernel $\{K(s,t),s,t\in T\}=\{u(s,t),s,t\in T\}$.
  In this case we refer to the permanental process  as an  associated  $\al$-permanental process. (It is associated with the  transient Markov process.) We use this terminology in what follows.
   
    We can use Theorem \ref{theo-9.2} to give    conditions for a permanental process to be unbounded in terms of the diagonals of the $M$-matrices  of  its finite dimensional distributions. Let $X=\{X_{ t},t\in T\}$\label{AA}, $T$ a countable set, be an $\al$-permanental process with  kernel $\{u(s,t),s,t\in T\}$.  Since, in  Theorem \ref{theo-9.2}, we require    that $A$ is an $M$-matrix, the  $\al$-permanental processes that we can consider must have a kernel  with the property  that for all $(t_{1},\ldots,t _{n})$ in $T$, the matrix with elements $\{u(t_{i},t_{j})\}_{i,j=1}^{n}$ is invertible and  its inverse  $A(t_{1},\ldots,t_{n})$  is a non-singular $M$-matrix.  This  is the case if (and only if)   $X$  is an associated $\al$-permanental process.  (This result is part of
  \cite[Theorem 13.1.2]{book}.     This  theorem it is stated for symmetric kernels but symmetry is not used in the proof. For the convenience of the reader,  in Section  \ref{subsec-7.1}, we repeat the proof of the portion of   \cite[Theorem 13.1.2]{book} that we use in this paper.)

  Suppose that $X$  is an associated $\al$-permanental process. Let $a_{i}(t_{1},\ldots,t_{n})$, $i=1,\ldots,n$, denote the diagonal elements of   $A(t_{1},\ldots,t_{n})$.
We use Theorem \ref{theo-9.2} in the following lemma which is proved in   Section \ref{proofs2}. 
It is a useful generalization of   Theorem \ref{theo-9.2} that enables us to only consider a fraction of the diagonal elements of $A$.

 \begin{lemma} \label{lem-1.2}
 Let $a^{*}_{i}(t_{1},\ldots,t_{n})$ denote   a non-decreasing rearrangement of  \newline $a_{i}(t_{1},\ldots,t_{n})$. For   any integer $p\ge 1$ let
 \begin{equation}
   \psi_{[n/p]}^{*}=\inf_{ (t_{1},\ldots,t _{n})\in T^{n}}a^{*}_{ [n/p]}(t_{1},\ldots,t_{n}).
   \end{equation}
Then  
\begin{equation}
      P\(\sup_{t\in T}  X_{t}\ge \la/  \psi_{ [n/p]}^{*}  \)\ge  P\(\max_{1\le i \le  [n/p]}  \xi^{(i)}_{ \al,1} \ge  \la \).\label{1.22e}
   \end{equation}
       \end{lemma}
       Therefore, if
       \begin{equation}
  \limsup_{n\to \ff}  P\(\max_{1\le i \le  [n/p]}  \xi^{(i)}_{ \al,1} \ge  \la_{n} \)=1\label{1.23w}
   \end{equation}
we have \begin{equation}
   P\(\sup_{t\in T}  X_{t}\ge \la_{n}/  \psi_{ [n/p]}^{*},\,\, i.o.  \)=1\label{n3.4n}.
   \end{equation}

  In Section \ref{proofs2} we show that (\ref{1.23w}) holds with $\la_{n}=\log n$. Therefore we can use   (\ref{n3.4n}) to     obtain   the following   sufficient  condition for   permanental processes to be unbounded.
 
 \begin{theorem}\label{theo-1.3} Let $\{X_{t}, t\in T\}$ be an   associated $\al$-permanental process.    If  
 \begin{equation}
\limsup_{n\to \ff}  \frac{\log n}{  \psi_{ [n/p]}^{*}}=\ff, \label{1.30ww}
   \end{equation}
 then $\sup_{t\in T}  X_{t}=\ff$ almost surely.     
 \end{theorem} 
 
The next  corollary is an immediate consequence of Theorem \ref{theo-1.3}.   
 
 \begin{corollary} \label{cor-1.2}
Let $d:=d_{s,t}$ be a    function on $T\times T$. Set 
   \be
 d_{n}^{ *} (t_{1},\ldots,t_{n})=\inf_{ 1\le i,j \le n,i\ne j}d _{t_{i},t_{j}}.\label{2.6q}
  \ee 
  Suppose that  a fraction  of the entries
\begin{equation}
   a_{i}(t_{1},\ldots,t_{n})\le   \frac{C'}{ (d_{n}^{ *})^{2}(t_{1}.\ldots,t_{n})},\label{1.34nn}
   \end{equation}
  for some constant $C'$. Then
   \begin{equation}
\limsup_{n\to \ff} \(\sup_{(t_{1},\ldots,t_{n})}  (d_{n}^{ *})^{2}(t_{1},\ldots,t_{n})\)\log n=\ff  ,
\label{sudxn}
 \end{equation}
imlies that  $\sup_{t\in T}  X_{t}=\ff$ almost surely. 
 \end{corollary}

  The condition in (\ref{1.34nn}) is not very useful because, in general one doesn't know the inverse of the matrices $\{u(t_{i },t_{j}\}_{i,j=1}^{n}$. In Lemma \ref{lem-5.2n} we give conditions on the  kernel $u(x,y)$ so that (\ref{1.34nn}) holds with the function 
\begin{equation}
   \si_{s,t}=\ (u(s,s)+u(t,t)-(u(s,t)+u(t,s))\ )^{1/2}\label{1.28-}
   \end{equation}
  replacing $d_{s,t}$. This enables us to obtain the following theorem:  
  
 \begin{theorem}\label{theo-5.1qq} Let $u$ be the potential density of a  transient Markov process
  in $R^{1}$ and assume that $u(s,s)$ is constant for all $|s|\leq \ep$,   for some $\ep>0$. Set 
\begin{equation}
   \si_{s,t}^{2} =2u(0,0)- u(s,t)- u(t,s)\label{levmet}
   \end{equation}
   and assume that 
        \begin{equation}
 |u(s,t)- u(t,s)|\le C\si_{s,t}^{2} , \qquad C<1,\label{4.20wwe}
 \ee
 for all $|s|,|t|\leq \ep$.    Then
    \begin{equation}
\limsup_{n\to \ff} \(\sup_{\stackrel{(t_{1},\ldots,t_{n})}{\forall |t_{i}|\leq \ep}}  (\si_{n}^{ *})^{2}(t_{1},\ldots,t_{n})\)\log n=\ff, 
\label{sudx}
 \end{equation}
  implies that the  $\al$-permanental  process with kernel  $u$  is unbounded  almost surely. 
 \end{theorem}
 
  It follows from Lemma \ref{lem-5.4} and that fact that $u(s,s)$ is constant for all $|s|\leq \ep$,   for some $\ep>0$, that (\ref{4.20wwe}) always holds for $C=1$.

 \medskip		In Theorem \ref{theo-5.1ee} we remove the hypothesis that $u(s,s)$ is constant for all $|s|$ sufficiently small. We don't consider this here because the result is not as  easy to state as  Theorem \ref{theo-5.1qq}.

 \medskip	   If $ u(s,t)$ is symmetric and positive definite it is the covariance of a Gaussian process.
Let  $\{\wt X_{t}, t\in R^{1}\}$ be a mean zero Gaussian process with covariance $  u(s,t) $. In this case
 \begin{equation}
   \si^{2}_{s,t}= E\(\wt X_{t}-\wt X_{s}\) ^{2}=u(s,s)+u(t,t)-2u(s,t).
   \end{equation}
(In particular this shows that $  \si_{s,t}$ is a metric on $R^{1}$.) 

Since  $\{\wt X_{t}, t\in R^{1}\}$ is a mean zero Gaussian process we can use  Slepian's Lemma to show that  (\ref{sudx}) implies that
 $\sup_{t\in R^{1}} \wt X_{t}=\ff$ almost surely.    This also follows from Theorem \ref{theo-5.1qq},  when $E\wt X_{t}^{2}$ is constant,  since in this case the left-hand side of (\ref{4.20wwe}) is equal to 0. (What Theorem \ref{theo-5.1qq} shows is that the 1/2-permanental process  $\sup_{t\in R^{1}} \wt X^{2}_{t}=\ff$ almost surely. Of course we also require that the inverse of $\{u(x_{j},x_{j})\}_{i,j=1}^{n}$ is an $M$-matrix for all $x_{i_{1}},\ldots,x_{i_{n}}\in R^{1}$.) 
 
Even when $ u(s,t)$ is not symmetric, $u(s,t)+u(t,s)$ is symmetric, and if it is also positive definite it is the covariance of a Gaussian process. In this case  we can still associate a permanental process with a Gaussian process. We plan to take this up in a subsequent paper.

 \medskip  We   can use Theorem \ref{theo-5.1qq} to study the boundedness of  permanental processes with kernels that are the potential densities of  transient L\'evy processes   in $R^{1}$. 	 Let $  Y=\{ Y_{t},t\in R_{+} \}$ be a L\'evy process and consider the transient   L\'evy  process  $\ov Y=\{\ov Y_{t},t\in R_{+} \}$    that is $Y$ killed at    $\xi_{1,1/\bb}$, an independent exponential time with mean $\bb>0 $.   If $u^{ \bb}(x,y)$ is the $\bb-$potential density of $Y$
 it is the zero potential of $\ov Y$ and thus is also the kernel of a permanental process.   In this example  $u^{ \bb}(x,y)=u^{ \bb}(0,y-x)=:u^{ \bb}(y-x )$.
  
  As we have mentioned above,  since
  $u^{\bb}(x,y)$ is the $0$-potential density of a  transient L\'evy process,  
 for every finite collection   $\,x_{ 
1},\ldots, x_{ n}\,\in\, R^{1}$, the $n\times n$ matrix
$U=  \{u(x_{ i},x_{ j} )\}_{ 1\leq i,j\leq n}$ is invertible and
its inverse is a non-singular $M$-matrix.   We use   Theorem \ref{theo-5.1qq} to find conditions under which  the $\al$-permanental  process with kernel $u^{\bb}$ is unbounded.

 We  write the characteristic function of    $ Y$ as 
 \begin{equation}
   E e^{i\la  Y_{t}}=e^{-t\psi(\la)}.
   \end{equation}
When $u^{ \bb}(y-x)$ is not symmetric,  $\psi(\la)$ is complex.  
 	 Set 
\begin{equation}
   \RR_{\bb}(\la)=    \RR e\(1/(\bb+\psi(\la))\)\quad\mbox{and}\quad    \II_{\bb}(\la)= \II m\(1/(\bb+\psi(\la))\).\label{def-RI}
   \end{equation}
   
  \begin{lemma} \label{lem-5.1}\cite[Lemma 5.2]{MRK}\label{lem-10.2} For $\bb>0$, assume that $\RR  _{\bb}(   \la)\in L^{1}(R_{+})$. Then the the $\bb$-potential density of $  X$ is
\begin{equation}
 u^{ \bb}(z)=  R_{\bb}(z)+  H_{\bb}(z) \quad\mbox{and}\quad  u ^{\bb}(-z)=R_{\bb}(z)-  H_{\bb}(z)\label{10.2ss},
 \end{equation}
 where
\begin{equation}
 R_{\bb}(z)=    \frac1{ \pi} \int_{0}^\ff \cos (\la z)\RR_{\bb}  (\la)\,d\la\label{10.3ss}
 \end{equation}
 and
\begin{equation}
 H_{\bb}(z)=  \frac1{ \pi} \int_{0}^\ff\sin (\la z)\II_{\bb}  (\la)\,d\la.\label{10.4ss}
 \end{equation}   
\end{lemma}

  As a special case of   (\ref{levmet}) we consider the metric  
  \bea
   \si _{\bb}({z })&=&\ (2u^{\bb}(0) - u^{\bb}(z )-u^{\bb}(-z) \ )^{1/2}\label{1.28qq}\\
   &=& \( \frac2{ \pi} \int_{0}^\ff (1-\cos (\la z))\RR_{\bb}  (\la)\,d\la\)^{1/2}\nn.
   \eea
 (Note that because   $\RR  _{\bb}(   \la)$ is positive and in $L^{1}(R_{+})$,  $ \Ga(x,y)  =R_{\bb}(y-x)$ is the covariance function of a stationary Gaussian process, say $\{G(z), z\in R ^{1}\}$. Therefore, $\si^{2}(z)=E(G(z)-G(0))^{2}$.)
 
 \medskip	 The following  condition for the $\al$-permanental  process with kernel $u^{\bb}$ to be unbounded is an immediate application of Theorem \ref{theo-5.1qq}.

  \begin{theorem} \label{theo-new}
   Suppose that  $\RR  _{\bb}(   \la)\in L^{1}(R_{+})$  and  
\begin{equation}
   |H_{\bb}(z)|\le C \si_{\bb}^{2}(z)\qquad\mbox{for some }\quad C<1/2\label{3.25nn}
   \end{equation}
and  all $|z|$ sufficiently small. Suppose, in addition, that  $\si_{\bb}^{2} (z) \ge f(|z|)$ for some increasing function $f$ for all $|z|$ sufficiently small.   Then  
\begin{equation}
   \limsup_{n\to \ff} f(1/n)\log n=\ff  \label{5.7qq}
   \end{equation}
implies that  the  $\al$-permanental  process with kernel $u^{\bb}$  is unbounded  almost surely. 
\end{theorem}

\begin{theorem} \label{theo-1.6}
Let $X=\{X(t),t\in R_{+}\}$ be  the $\al$-permanental process with a kernel  that is the the $\bb$ potential density of a L\'evy process with L\'evy measure
 \begin{equation}
   \nu(dx)=\(x^{-2}  g(1/|x|)(pI_{x>0}+qI_{x<0})\)dx \quad p,q> 0, \quad p+q=1,\label{1.30s}
   \end{equation}
in which $g$ is a positive,   quasi-monotonic slowly varying function at infinity.  Suppose $p\ne q$ and
\begin{equation}
  \lim_{n\to\ff} \int_{1}^{n}    \frac{g(s)}s\,ds= \ff.\label{1.30}
   \end{equation}Then $X$ is unbounded almost surely if  
    \bea 
   \int_{1}^{n}    \frac{g(s)}s\,ds=o(\log n)\label{6.19s},
 \eea
   as $n\to\ff$. 
   
  If $p=q$ and  
   \begin{equation}
    \lim_{n\to\ff} \int_{n}^{\ff}\frac{1}{sg(s)}\,ds=0,
   \end{equation}
   then $X$ is unbounded almost surely if 
   \begin{equation}
    \( \int_{n}^{\ff}\frac{1}{sg(s)}\,ds\)^{-1}=o(\log n).\label{1.33}
   \end{equation} 
 \end{theorem}

 It is interesting to note that the $\bb$ potential density determined by (\ref{1.30s}) has the property that 
for $z>0$  
\bea
   u^{\bb}(z)&\sim&  u^{\bb}(0)-\frac{\si^{2}(z)}2\(1-   |p-q| \)\label{1.31}\\
      u^{\bb}(-z)&\sim&  u^{\bb}(0)- \frac{\si^{2}(z)}2\(1+  |p-q|  \).\nn
      \eea
  as $z\to 0$.   Here we write $f\sim g$ as $z\to 0$ if $\lim_{ z\to 0}f(z)/g(z)=1$, with a similar meaning for 
$f\sim g$ as $z\to \ff$.   The derivation of  (\ref{1.31}) is given in Section \ref{sec-6} following the proof of Theorem 
\ref{theo-1.6}.

\begin{example}\label{ex-1.1} {\rm 

    We consider Barlow's example \cite[page 1393]{Barlow},   slightly modified,  of a L\'evy process  with  L\'evy measure given by (\ref{1.30s}) with $g(y)$ replaced by $   g_{\ga\de}(y)$ where
 \begin{equation}
   g_{\ga\de}(y)=(\log y)^{\ga}(\log \log y)^{\de}1_{\{y>\ep\}},
   \end{equation}with $\ga>-1$.   
  Let $Y_{\ga\de}$ be the  L\'evy process determined by this L\'evy measure  and denote its $\bb$ potential density by $u^{\bb}$. It follows from (\ref{6.19s})  that when $p\neq q$  the permanental process   with kernel $u^{\bb}$ is unbounded    if    $\ga< 0$ or $\ga=0$ and $\de<0$.

 \medskip	  When $p=q$, 
   \begin{equation}
    \( \int_{n}^{\ff}\frac{1}{s(\log s)^{\ga}(\log \log s)^{\de}  }\,ds\)^{-1}\sim C (\log n)^{\ga-1}(\log \log n)^{\de}.\label{1.33q}
   \end{equation}  
and we now require that $\ga>1$.  In this case the permanental process with kernel $u^{\bb}$ is unbounded    if  $\ga<2$  or $\ga=2$ and $\de<0$.
}\end{example}

  Let $u^{\bb}(s,t)=u^{\bb}(t-s)$ be the $\bb$-potential of a L\'evy process. Using Barlow's \cite{ Barlow} necessary and sufficient condition for the boundedness of local times of L\'evy processes and an isomorphism theorem of Eisenbaum and Kaspi \cite{EKC}, that relates local times and permanental processes we can show that the associated $\al$-permanental process is unbounded almost surely if the Gaussian process with covariance $\ga(s,t)=u^{\bb}(s-t)+ u^{\bb}( t-s))$   is unbounded almost surely.  (See the comment following Lemma \ref{lem-5.1}.) For the processes considered in Example \ref{ex-1.1} this occurs if and only if
\begin{equation}
   \int_{1}^{\ff}\frac{\(\int_{\la}^{\ff}\RR _{\bb}(u)\,du\)^{1/2}}{ \la(\log\la)^{1/2}}\,d\la=\ff.
   \end{equation}
  Consequently, when $p\ne q$ the  the permanental process   with kernel $u^{\bb}$ in Example \ref{ex-1.1} is unbounded  almost surely   if $\ga<0$ or $\ga=0$ and $\de\le 2$ and bounded almost surely when $\ga=0$ and $\de>2$. When $p=q$ it is unbounded  almost surely   if $\ga<2$ or $\ga=2$ and $\de\le 2$ and bounded almost surely when $\ga=2$ and $\de>2$. This gives a little more than we obtain in Example \ref{ex-1.1}. 
 
\medskip	  Even though the results in  Theorems \ref{theo-new} and \ref{theo-1.6} are not best possible, the theorems are interesting for at least two reasons. The first is that their proofs are much simpler than the proof in \cite{Barlow}. The second is that the proofs involving 
\cite{ Barlow} and \cite{EKC}  are indirect and give  no insight into why permanental processes have sample path properties similar to the squares of  Gaussian processes. Our proofs of  Theorems \ref{theo-new} and \ref{theo-1.6} are classical and relatively simple and show  that permanental processes have sample path properties similar to the squares of  Gaussian processes because the Permanental 
Inequality is a generalization, in many respects, of the Sudakov Inequality.

\medskip	 With some restrictions and a simplification, and slight weakening, of  (\ref{3.25nn}) we get a Corollary of Theorem \ref{theo-new} that is easier to use  and imposes weaker conditions on the behavior of  $ |\II  _{\bb}(   \la) |$ and $ \RR  _{\bb}(   \la) $ as 
$\la\to\ff $.
     \begin{corollary} \label{cor-1.4}
 \label{theo-5.1}Suppose that   $\RR  _{\bb}(   \la)\in L^{1}(R_{+})$ and that  $ |\II  _{\bb}(   \la) |$ and $ \RR  _{\bb}(   \la) $ are asymptotic to  non-increasing functions as $\la\to\ff$ and  
         \begin{equation}
 |   z|\int_{0}^ {\pi/|z|}  \la \, | \II  _{\bb}(   \la)| \,d\la   \le \frac{C}{2} \int_{\pi/(2|z|)}^\ff  \RR_{ \bb}( \la) \,d\la \label{5.11}
   \end{equation}
 for some $C<1$, for all $|z|$ sufficiently small. Then  
  \begin{equation}
 \limsup_{n\to \ff} \(  \int_{n}^\ff  \RR_{ \bb}( \la) \,d\la\)\log n=\ff  \label{5.7q}
   \end{equation}
implies that  the $\al$-permanental  process with kernel $u^{\bb}$ is unbounded  almost surely. 
 \end{corollary}

 \medskip	  The proof of Theorem \ref{theo-cup} is given in Section \ref{proofs} and that of Theorem  \ref{theo-1.3} in Section \ref{proofs2}.  In Section  \ref{sec-10} we examine the implications of (\ref{10.14}), the Permanental Inequality  and explain why we refer to it as a Sudakov type inequality.    Theorem \ref{theo-5.1qq} is proved in Section \ref{sec-5}. In Section \ref{sec-6} we prove Theorem \ref{theo-new},  Corollary \ref{cor-1.4} and fill in the details for Example \ref{ex-1.1}.

\section{Representation of permanental processes }\label{proofs}

\medskip	For any $n\times n$ matrix $M$   we define the  $\al$-permanent  
\begin{equation}
|M|_{\al}= \begin{vmatrix}
 m_{1,1} & \cdots&  m_{1,n}\\ 
\cdots& &\cdots \\
 m_{n,1} &\cdots& m_{n,n}
\end{vmatrix}_{\al}=\sum_{\pi}\al^{c(\pi)} m_{1,\pi(1)} m_{2,\pi(1)}\cdots  m_{n,\pi(n)}.\label{vj.71}
\end{equation}
Here the sum runs over all permutations $\pi$ on $[1,n]$ and $c(\pi)$ is the number of cycles in $\pi$.
We make the trivial observation that if all entries of $M$ are non-negative, then $|M|_{\al}\geq 0$. 

 We use boldface, such as  ${\bf x }$, to denote vectors. Let   ${\bold k}=(k_{1},\ldots, k_{n})\in \mathbb{N}^{n}$ and $ |\bold k|=\sum_{l=1}^{n}k_{l}$. For  $1\leq p\leq  |\bold k|$,   set    $ i_{p}=j$, where 
 \begin{equation}
\sum_{l=1}^{j-1}k_{l}< p\leq \sum_{l=1}^{j}k_{l}.\label{vj.70v}
 \end{equation}
For any $n\times n$ matrix  $C=\{ c_{ i,j}\}_{ 1\leq i,j\leq n}$  we define  
  \begin{equation}
C(\bold k)= \begin{bmatrix}
c_{ i_{1}, i_{1}} & c_{ i_{1}, i_{2}}&\cdots &c_{i_{1}, i_{|\bold k|}}\\ 
c_{ i_{2}, i_{1}} & c_{ i_{2}, i_{2}}&\cdots &c_{i_{2}, i_{|\bold k|}}\\ 
\cdots& &\cdots \\
c_{ i_{|\bold k|}, i_{1}} & c_{ i_{|\bold k|}, i_{2}}&\cdots &c_{i_{|\bold k|}, i_{|\bold k|}}
\end{bmatrix},\label{vj.71va}
\end{equation}
 and $C(\bold 0)=1$.	For example,   if $n=3$ and  ${\bold k}=(0,2, 3)$ then $|\bold k|=5$ and $i_{1}=i_{2}=2$ and $i_{3}=i_{4}=i_{5}=3$.  
   \begin{equation}
C(0,2, 3)= \begin{bmatrix}
c_{ 2, 2} & c_{2, 2}&c_{2, 3}&c_{2, 3}&c_{2, 3}\\ 
c_{ 2, 2} & c_{2, 2}&c_{2, 3}&c_{2, 3}&c_{2, 3}\\  
c_{ 3, 2} & c_{3, 2}&c_{3, 3}&c_{3, 3}&c_{3, 3}\\ 
c_{ 3, 2} & c_{3, 2}&c_{3, 3}&c_{3, 3}&c_{3, 3}\\ 
c_{ 3, 2} & c_{3, 2}&c_{3, 3}&c_{3, 3}&c_{3, 3}
\end{bmatrix}.\label{vj.71vs}
\end{equation}

\begin{lemma} \label{lem-1.1}
Let $A$ be an $n\times n$ nonsingular  M-matrix with diagonal entries $a_{1},\ldots, a_{n}$ and   $S$ be an $n\times n$ diagonal matrix with entries $(s_{1},\ldots, s_{n})$ and set 
 \begin{equation}
A=D-B,\label{vj.73}
\end{equation}
where $D=\mbox{diag }(a_{1},\ldots, a_{n})$   and  all the elements of $B$ are non-negative, (so that   all the diagonal elements of $B$ are equal to zero). 
 Then
\bea 
 \frac{|A|^{\al}}{|A+S|  ^{ \al}}
 &=& |A|^{\al}\sum_{{\bold k}=(k_{1},\ldots, k_{n}) }  {|B(\bold k)|_{\al}\over k_{1}! \cdots k_{n}!}\,\frac{1}{ (a_{1}+s_{1})^{\al+k_{1}}\cdots(a_{n}+s_{n})^{\al+k_{n}}}\nn\\
  &=& \frac{|A|^{\al}}{\prod_{i=1}^{n}a^{\al}_{i}}\sum_{{\bold k}=(k_{1},\ldots, k_{n}) }  {|B(\bold k)|_{\al}\over  \prod_{i=1}^{n}a_{i}^{k_{i}} k_{i}! }\prod_{i=1}^{n}\,\(\frac{a_{i}}{  a_{i}+ s_{i}  }\)^{\al+k_{i} }  ,\label{9.9}
   \eea
where   the sum is over all ${\bold k}=(k_{1},\ldots, k_{n})\in \mathbb{N}^{n}$.   (The series converges for all $s_{1},\ldots,s_{n}\in R_{+}^{n}$ for all $\al>0$.)
  \end{lemma}
  
\Proof  For $B$ as given in (\ref{vj.73}) consider
\begin{equation}
 \HH(z_{1},\ldots,z_{n}) = |I-ZB|^{-\al}= \sum_{{\bold k}=(k_{1},\ldots, k_{n}) }\(\prod_{i=1}^{n}\frac{z_{i}^{k_{i}}}{k_{i}!}\)|B(\bold k)|_{\al},\label{vj.75s}
   \end{equation}
where $Z$ is a diagonal matrix with entries $z_{1},\ldots,z_{n}$   and the second equality is   given in \cite[(6)]{VJ}. By \cite[Theorem, page 120]{VJA}, the   series  (\ref{vj.75s}) converges for $(z_{1},\ldots,z_{n})$, when the modulus of the maximum eigenvalue of $ZB$ is less than $1$. 

 We write 
\begin{equation}
|A+S| =|(D+S)-B| =|(D+S)| |I-(D+S)^{-1}B|,\label{vj.74}
\end{equation}
so that 
\be  |A+S|^{-\al}
\label{vj.75t}  = \sum_{{\bold k}=(k_{1},\ldots, k_{n}) }  {|B(\bold k)|_{\al}\over k_{1}! \cdots k_{n}!}\,\frac{1}{ (a_{1}+s_{1})^{\al+k_{1}}\cdots(a_{n}+s_{n})^{\al+k_{n}}}.
\ee
By the statements in the first paragraph of this proof, this series converges when the modulus of the maximum eigenvalue of $(D+S)^{-1}B$ is less than $1$.   

We complete the proof by referring to several results in the  valuable book \cite{BP}. Note that the definition of $M$-matrix on \cite[pg. 133]{BP}  is different from the one that we give. However, it follows by   \cite[$N_{38}$, pg. 137]{BP} that they are equivalent.
We now write $A+S=D+S-B$, to see by  \cite[Chapter 7, Theorem 5.2]{BP}, that the maximum eigenvalue of $(D+S)^{-1}B$ is less than 1 if and only if $(A+S)^{-1}B\geq 0$. Since $A$ is a non-singular $M$-matrix,   we have  $B\ge 0$ and by  \cite[Chapter 6, Theorem 2.4]{BP}   $(A+S)^{-1}\ge 0$ as well. This completes the proof of this lemma.\qed

    In the next theorem we   give an explicit description of random variables with  Laplace transforms given in (\ref{9.9}).
    
      \begin{theorem} \label{theo-9.1}Let  $A$ be    an $n\times n$ non-singular $M$-matrix as in Lemma \ref{lem-1.1}.  Set    $Z=(Z_{1},\ldots,Z_{n})$ with
 \begin{equation}
  P\(Z=(k_{1},\ldots, k_{n})\)=\frac{|A|^{\al}}{\prod_{i=1}^{n}a^{\al}_{i}}  {|B(\bold k)|_{\al}\over  \prod_{i=1}^{n}a_{i}^{k_{i}} k_{i}! },\label{9.11w}
  \end{equation} 
  and $X=(X_{1},\ldots,X_{n})$ with  
 \bea
  X&=& \(\xi^{(Z,1)}_{ \al+Z_{1},a_{1}},\ldots,\xi^{(Z,n)}_{\al+Z_{n},a_{n}}\)\label{9.10}\\
&=&  \sum_{{\bold k}=(k_{1},\ldots, k_{n}) } 1_{   k_{1},\ldots, k_{n}  } (Z)\( \xi^{({\bf k},1)}_{ \al+k_{1},a_{1}},\ldots,\xi^{({\bf k},n)}_{\al+k_{n},a_{n}}\),\nn
  \eea  
  where $Z$ and all the gamma distributed random variables,  $\xi^{({\bf k},i)}_{\cd,\cd}$, ${\bf k}\in \mathbb{N}^{n}, i\in{1,\ldots,n}$    are independent  and   $\{a_{i}\}_{i=1}^{n}$ are the diagonal elements of $A$.  
Then   
\begin{equation}
   E\(e^{-\sum_{i=1}^{k}s_{i}X_{i}}\) 
 =   \frac{|A|^{\al}}{ |A+S |^{ \al}}. \label{9.15} 
 \end{equation}
  \end{theorem}

   \noindent {\bf Proof  of   Theorem \ref{theo-9.1}} Taking $S=0$ in (\ref{9.9}) we see that
\begin{equation}
 \frac{|A|^{\al}}{\prod_{i=1}^{n}a^{\al}_{i}}\sum_{{\bold k}=(k_{1},\ldots, k_{n}) }  {|B(\bold k)|_{\al}\over  \prod_{i=1}^{n}a_{i}^{k_{i}} k_{i}! }=1.\label{9.9g}
\end{equation}
Therefore we can define an  $\mathbb{N}^{n}$ valued random variable $Z=(Z_{1},\ldots,Z_{n})$ with
 \begin{equation}
  P\(Z=(k_{1},\ldots, k_{n})\)=\frac{|A|^{\al}}{\prod_{i=1}^{n}a^{\al}_{i}}  {|B(\bold k)|_{\al}\over  \prod_{i=1}^{n}a_{i}^{k_{i}} k_{i}! }.\label{9.11ww}
  \end{equation}
  We write   (\ref{9.9}) in the form
 \begin{equation}
   \frac{|A|^{\al}}{ |A+S |^{ \al}}=\sum_{{\bold k}=(k_{1},\ldots, k_{n}) } P\(Z=(k_{1},\ldots, k_{j})\)\prod_{i=1}^{n} \,\(\frac{a_{i}}{  a_{i}+ s_{i}  }\)^{\al+k_{i} } .\label{vj.58}
 \end{equation}
This is the   Laplace transform of the  $R^{n}_{+}$ valued  random variable  \bea
  X&=& \(\xi^{({\bf Z },1)}_{ \al+Z_{1},a_{1}},\ldots,\xi^{({\bf Z },n)}_{\al+Z_{n},a_{n}}\)\label{9.10q}\\
&=&  \sum_{{\bold k}=(k_{1},\ldots, k_{n}) } I_{   k_{1},\ldots, k_{n}  } (Z)\(\xi^{({\bf k },1)}_{ \al+k_{1},a_{1}},\ldots,\xi^{({\bf k },n)}_{\al+k_{n},a_{n}}\),\nn
  \eea
  where all the random variables   are  independent.
 \qed

  \medskip	 \noindent {\bf Proof of  Theorem \ref{theo-cup}  }  Theorem \ref{theo-cup}  follows  from (\ref{9.10})  and   the facts that
  \be
 \xi_{\al+k_{i},a_{i}}\stl\xi_{ \al,a_{i}}+\xi_{k_{i},a_{i}},\label{1.15a}
 \ee
 and
\be
 \xi_{\al,a_{i}}\stl a_{i}^{-1}\xi_{ \al,1},\label{1.15ar}
 \ee
 which allow us to write  
   \bea
   X&=& \(\xi^{({\bf Z },1)}_{ \al+Z_{1},a_{1}},\ldots,\xi^{({\bf Z },n)}_{\al+Z_{n},a_{n}}\)\label{1.15}\\
   &\stl&  \(a_{1}^{-1}\xi^{(1)}_{ \al,1},\ldots,a_{n}^{-1}\xi^{(n)}_{ \al,1}\)+ 
 \(\xi^{({\bf Z },1)}_{   Z_{1},a_{1}},\ldots,\xi^{({\bf Z },n)}_{ Z_{n},a_{n}}\),\nn
   \eea
   where $\xi^{(i)}_{ \al,1}$ are i.i.d. copies of $\xi _{ \al,1}$ and we set $\(\xi^{\cd}_{   0,a_{1}},\ldots,\xi^{\cd}_{ 0,a_{n}}\)=0$.  \qed

We get the following immediate corollary of Theorem \ref{theo-9.1}
  
\begin{corollary}\label{cor-infdiv} Let $A$ be   an $n\times n$ non-singular $M$-matrix. Then for each $\al>0$, (\ref{2.1}) defines an $n$-dimensional infinitely divisible random variable. \end{corollary}

Actually Eisenbaum and Kaspi \cite[Lemma 4.2]{EK} show that the condition in Corollary \ref{cor-infdiv} is both necessary and sufficient. They do this by extending the proof of this result by Bapat, Griffiths and Milne in the case when $K$ is symmetric,   (see \cite{EK} for references), to the case when $K$ is not symmetric. The proof of sufficiency in Corollary \ref{cor-infdiv} is completely different from their proof.

  \medskip	
  It follows from (\ref{9.11w}) and  (\ref{9.10}) that  for measurable functions  $f$ on  $R^{n}_{+}$, 
   \bea
 && E(f(X))\label{9.12}\\
  &&\qquad=\frac{|A|^{\al}}{\prod_{i=1}^{n}a^{\al}_{i}}\sum_{{\bold k}=(k_{1},\ldots, k_{n}) }  {|B(\bold k)|_{\al}\over  \prod_{i=1}^{n}a_{i}^{k_{i}} k_{i}! }E\( f\(\xi^{({\bf k},1)}_{ \al+k_{1},a_{1}},\ldots,\xi^{({\bf k},n)}_{\al+k_{n},a_{n}}  \)\).\nn
  \eea

\medskip	 Obviously  (\ref{9.12}) gives us more than (\ref{10.14}). Even though it is difficult to compute $B(\bold k)$ for all ${\bold k}$ it is not difficult to obtain it for some ${\bold k}$ and to improve   (\ref{10.14}).

    \medskip	  All the results in this paper follow from the representation in Lemma \ref{lem-1.1}. A different form of this representation, under different hypotheses, is given in \cite{MM}. It seems to be more useful than Lemma \ref{lem-1.1} in obtaining explicit probability density functions of low dimensional multivariate gamma distributions. Lemma \ref{lem-1.1} is more useful in describing multivariate gamma distributions in high dimensions. (Multivariate gamma random variables and $\al$-permanental random variables are synonyms.)

\section{Proof of Theorem \ref{theo-1.3}}\label{proofs2}
    
   The next three lemmas are used in the proof of Theorem \ref{theo-1.3}. 
    
  \begin{lemma} \label{lem-2.2qq}   For  $ \la>2(u-1) \vee 0 $   and all $u,v>0$
 \begin{equation}
  P\(\xi_{u,v}\ge\la/v\)   \le \frac{2  \la ^{u-1}e^{- \la}}{\Gamma(u)}\label{4.1e}.
   \end{equation}
   and for  $ \la\ge2 $   and all $u,v>0$
    \begin{equation}
 \frac{2  \la ^{u-1}e^{- \la}}{3\Gamma(u)} \le P\(\xi_{u,v}\ge\la/v\)   \label{4.1ee}.
   \end{equation}
    \end{lemma}
 
   \Proof    Using the fact that $P\(\xi_{u,v}\ge\la/v\) =P\(\xi_{u,1}\ge\la \)$ it suffices  to get the bounds in (\ref{4.1e}) for $P\(\xi_{u,1}\ge\la \)$.
 By an integration by parts 
  \begin{equation}
   \int _{\la}^{\ff}    x^{u-1} e^{-  x}  \,dx= \la^{u-1}e^{-\la}+(u-1) \int _{\la}^{\ff}    x^{u-2} e^{-  x}  \,dx\label{4.2j1}
  \end{equation}
  The upper bound in (\ref{4.1e}) follows immediately if $u\leq 1$.
If $u>1$ and    $\la>2(u-1) $ 
  \begin{eqnarray}
  (u-1) \int _{\la}^{\ff}    x^{u-2} e^{-  x}  \,dx &\leq & { \la\over 2}\int _{\la}^{\ff}    x^{u-2} e^{-  x}  \,dx 
  \label{4.2j3q}\\
&\leq & { 1\over 2}\int _{\la}^{\ff}    x^{u-1} e^{-  x}  \,dx. \nn
  \end{eqnarray}
 Using this in   (\ref{4.2j1}) we see that
\begin{equation}
   \int _{\la}^{\ff}    x^{u-1} e^{-  x}  \,dx\leq 2\la^{u-1}e^{-\la}.\label{4.2j4}
\end{equation}
This gives the upper bound in (\ref{4.1e}).
  
To obtain the lower bound we first note that   for $u\geq 1$ it follows from (\ref{4.2j1}) that for any $\la>0$ we have
    \begin{equation}
   \int _{\la}^{\ff}    x^{u-1} e^{-  x}  \,dx\geq \la^{u-1}e^{-\la}.\label{4.2j2}
  \end{equation} 
  Similarly,  for $u<1$,     we use  (\ref{4.2j1}) to see that for any $\la>0$  
     \begin{equation}
   \int _{\la}^{\ff}    x^{u-1} e^{-  x}  \,dx = \la^{u-1}e^{-\la}-(1-u) \int _{\la}^{\ff}    x^{u-2} e^{-  x}  \,dx,\label{4.2j5}
  \end{equation}
 and since, for     $\la>2(1-u) $ 
  \begin{eqnarray}
  (1-u) \int _{\la}^{\ff}    x^{u-2} e^{-  x}  \,dx &\leq & { \la\over 2}\int _{\la}^{\ff}    x^{u-2} e^{-  x}  \,dx 
  \label{4.2j3}\\
&\leq & { 1\over 2}\int _{\la}^{\ff}    x^{u-1} e^{-  x}  \,dx, \nn
  \end{eqnarray}
we get
 \begin{equation}
   \int _{\la}^{\ff}    x^{u-1} e^{-  x}  \,dx\geq { 2\over 3}\la^{u-1}e^{-\la}.\label{4.2j6}
\end{equation}
Combining (\ref{4.2j2}) and (\ref{4.2j6}) we get the lower bound in (\ref{4.1e}).\qed

\begin{lemma} \label{lem-4.2}Let $\{\xi^{(i)}_{u,v}\}_{i=1}^{n}$ be independent. Then  for all   $\ep,q>0$,    $n\ge 10$ and $ (n^{\ep}/(q \,\Gamma(u)  \log n) )\ge 3/2$,
 \begin{equation}
   P\(\max_{1\le i\le n} \xi^{(i)}_{u,v  }\ge  \frac{ (1-\vep)\log n}{v }\)\ge 1-e^{-q}.\label{4.5}
   \end{equation}
 \end{lemma}
 
  \Proof  
We have  
  \bea
   P\(\max_{1\le i\le n} \xi^{(i)}_{u,v  } > \frac{(1-\ep)\log n}{v  }\)\   &=& 1-P\(\max_{1\le i\le n} \xi^{(i)}_{u,v  } \leq  \frac{(1-\ep)\log n}{v }\)\label{4.8}\\& =&1-\prod_{i=1}^{n }\(1- P\( \xi^{(i)}_{u,v  }> \frac{(1-\ep)\log n}{v  }\) \)\nn.
   \eea
  By  (\ref{4.1ee}),    for $n^{\ep}/(q  \Gamma(u) \log n)\ge 3/2$, 
\begin{equation}
   P\( \xi^{(i)}_{u,v  }> \frac{(1-\ep)\log n}{v  }\)\ge 
  \frac{2e^{-(1-\ep)\log n}} {3\Gamma(u)(1-\ep)\log n}\ge \frac{q}{n}.
   \end{equation}
 Using this and (\ref{4.8}) we see that   
    \begin{equation}
P\(\max_{1\le i\le n} \xi^{(i)}_{u,v  } > \frac{(1-\ep)\log n}{v  }\)\ge  1- \(1-  \frac{q}{n}  \)^{n }>1-e^{-q} \nn.
   \end{equation}
   \qed

The next lemma follows immediately from (\ref{10.14}). It is useful because in applying the Permanental Inequality sometimes  we don't want to consider all the diagonal elements of the non-singular $M$-matrix $A$. We use it in the proof of Lemma \ref{lem-1.2}

\medskip	For a sequence $\{v_{i}\}_{i=1}^{k}$ we define  $\{v^{*}_{i}\}_{i=1}^{k}$ to be the non-decreasing rearrangement of $\{v_{i}\}_{i=1}^{k}$.    

  	  \begin{lemma}  \label{lem-4.3}Let $X=(X_{1}.\ldots, X_{n})$ be an $R^{n}$ valued random variable defined by (\ref{2.1}) with an $n\times n$ non-singular $M$-matrix $A$  with diagonal elements $a_{i}$, $1\le n$.  Then for all   $p\ge 1$,  
     \be
   P\(\max_{1\le i\le n} X_{ i}  \ge \la \)\ge  P\((a^{*}_{[n/p] })^{-1} \max_{1\le i\le [n/p] }  \xi^{(i)}  _{ \al,1} \ge \la \),\label{1.11aa}   \ee
 where    $ \{   \xi ^{(i)} _{ \al,1}$, $1\le i\le [n/p]\}$ are independent.    \end{lemma}
  
  \Proof   
Using  (\ref{10.14}) we see that 
      \bea
   P\(\max_{1\le i\le n} X_{ i}  \ge \la \)&\ge &P\(\max_{1\le i\le n } a^{-1} _{i } \xi^{(i)}  _{ \al,1} \ge \la \)\label{4.10qj} \\
   &=&   P\(\max_{1\le i\le n }  (a^{\ast})^{-1} _{i } \xi^{(i)}  _{ \al,1} \ge\la \)\nn\\&\ge &P\(\max_{1\le i\le [n/p] } (a^{*} _{i })^{-1} \xi ^{(i)} _{ \al,1} \ge \la \)\nn \\
   &\ge &  P\((a^{*} _{[n/p] })^{-1}\max_{1\le i\le [n/p]}  \xi ^{(i)} _{ \al,1} \ge\la \) \nn.
   \eea
  \qed 
  
  \noindent {\bf Proof  of   Lemma  \ref{lem-1.2} }  By Lemma \ref{lem-4.3}, for any sequence $t_{1},\ldots,t_{n}\in T$,   
    \bea 
   P\(\sup_{t\in T}X_{t}\ge\la\)&\ge &P\(\max _{1\le i\le n}X_{t_{i}}\ge\la\)\label{3.15}\\
   &\ge& P\((a^{*}_{{[n/p]}}(t_{1},\ldots,t_{n}))^{-1}\max_{1\le i\le [n/p]}  \xi ^{(i)}  _{\al,1} \ge\la \)\nn\\
   &\ge& P\(\max_{1\le i\le [n/p]}  \xi ^{(i)}  _{\al,1} \ge  a^{*}_{{[n/p]}}(t_{1},\ldots,t_{n})  \la \)\nn.
  \eea
 Therefore, by continuity  of the cumulative distribution function \newline $ P\(\max_{1\le i\le [n/p]}  \xi ^{(i)}  _{\al,1} \le s  \)$, we have  
 \begin{equation}
    P\(\sup_{t\in T}X_{t}\ge\la\)\ge  P\(\max_{1\le i\le [n/p]}  \xi ^{(i)}  _{\al,1} \ge \inf_{(t_{1},\ldots,t_{n})\in T^{n}}   a^{*}_{{[n/p]}}(t_{1},\ldots,t_{n})  \la \), 
   \end{equation}
which is (\ref{1.22e}).\qed
  
 \noindent {\bf Proof  of   Theorem \ref{theo-1.3} }   By 
(\ref{1.23w})  all we need to do is to show that    
\begin{equation}
  \lim_{n\to\ff}  P\(\max_{1\le i \le  [n/p]}  \xi^{(i)}_{ \al,1} \ge  \log n \)=1.
   \end{equation}
This follows immediately from Lemma \ref{lem-4.2}.\qed

 	\section{Permanental Inequality}\label{sec-10} 
	
	We examine the implications of (\ref{10.14}) and explain why we refer to it as a Sudakov type inequality.

  	 It follows from  (\ref{10.14}), the Permanental Inequality,  that
     \begin{equation}
   E\( \max_{1\le i\le n} (2X_{ i})^{1/2}  \)\ge E\(\max_{1\le i\le n } (2 \xi ^{(i)}_{ \al,1}/a _{i})^{1/2}  \). \label{2.2s}
   \end{equation}
  If $K=A^{-1}$ is symmetric and positive definite and $\al=1/2$,  then $X=(\eta_{1}^{2}/2 ,\ldots, \eta_{n}^{2}/2)$, where $(\eta_{1},\ldots,\eta_{n})$
is a Gaussian vector with covariance $\{K_{i,j}\}$. In this case  by (\ref{2.2s}) 
\be 
   E\( \max_{1\le i\le n} |\eta_{i} |  \ge \la \)\ge \sqrt2  E\(\max_{1\le i\le n } ( \xi ^{(i)}_{ \al,1}/a _{i})^{1/2}  \). \label{2.3}
      \ee   
  
   Note that $ \(   \xi ^{(i)}_{ \al,1}\)^{1/2}$, $1\le i\le n$, are the absolute values of a sequence of independent normal random variable with variance   $1/2$. Therefore we can rewrite (\ref{2.3}) as   
      \begin{equation}
   E\( \max_{1\le i\le n} | \eta_{i}|   \)\ge\sqrt2 E\(\max_{1\le i\le n }  |\zeta_{i}| /\sqrt{2a_{i}}      \)\ge  \( {\max_{1\le i\le n }\sqrt{a_{i}}}\)^{-1}E\(\max_{1\le i\le n }  |\zeta_{i}|       \),\label{2.5s}
   \end{equation}
where $\ze_{i}$, $1\le i\le n$, are independent normal random variables with mean zero and variance 1. This is what we get from the Permanental Inequality for a mean zero normal random vector $(\eta_{1},\ldots,\eta_{n})$ with covariance matrix $K$.

 	By Fernique's comparison principle \cite[Lemma 5.5.3]{book}
          \begin{equation} 
   E\( \max_{1\le i\le n}  \eta_{i}  \)\ge E\(\max_{1\le i\le n }  \rho_{i}   \),\label{2.4s}
   \end{equation}
   where $(\rho_{1},\ldots,\rho_{n})$ is a mean zero Gaussian random variable
   satisfying
   \begin{equation}
   E(\rho_{i}-\rho_{j})^{2}\le E(\eta_{i}-\eta_{j})^{2} = K_{i,i} +K_{j,j}-2K_{i,j} .\label{prob5}
   \end{equation}
   This can be achieved when
  $   \rho_{i  } $, $1\le i\le n$, are   independent normal random variable with variance $(\si^{*}_{n})^{2}/2$ where 
  \be
  (\si_{n}^{ *})^{2}=\inf_{ 1\le i,j \le n,i\ne j} K_{i,i} +K_{j,j}-2K_{i,j}.\label{2.6}
  \ee 
  With this choice of   $   \rho_{i  } $, $1\le i\le n$ we get
       \begin{equation} 
   E\( \max_{1\le i\le n}  \eta_{i}  \)\ge \frac{ \si_{n}^{ *} }{\sqrt 2}E\(\max_{1\le i\le n }  \ze_{i}   \).\label{2.4ss}
   \end{equation}
   This inequality is essentially Sudakov's Inequality.

If we  ignore  the presence or absence of absolute values we see that  if
  \begin{equation}
   \max_{1\le i\le n} a_{i}\le \frac{2}{(\si_{n}^{ *} )^{2}},\label{2.6aaa}
   \end{equation}   
  then (\ref{2.5s}), which   follows  from the Permanental Inequality, gives a stronger lower bound for $  E\( \max_{1\le i\le n}  \eta_{i}    \)$ than   (\ref{2.4ss}), which is what we get using the Sudakov Inequality.  
  In Lemma \ref{lem-5.2n} we show that (\ref{2.6aaa}) holds when the matrix $K$ is   symmetric and constant on the diagonals. 

\begin{remark} {\rm The Sudakov Inequality is very useful in giving necessary conditions for a Gaussian process to be bounded, but it can be a very weak  lower  bound for many Gaussian random variables. We point this  out because the Permanental Inequality has the same limitations. Evaluating the right-hand side of (\ref{2.4ss}) we get
     \begin{equation} 
   E\( \max_{1\le i\le n}  \eta_{i}  \)\ge C \si_{n}^{ *}   (\log n)^{1/2}.\label{3.10}
   \end{equation}
for some constant $C>0$, for all $n$ sufficiently large. If we take the limit as $n\to\ff$, as we do when considering whether a Gaussian process is bounded, this is only useful when
\begin{equation}
   \limsup_{n\to\ff}\si_{n}^{ *}   (\log n)^{1/2}>0.
   \end{equation}
   Let $\{B(t),t\in [0,1]\}$ be Brownian motion and consider $(B(1/n),B(2/n),\ldots,\newline	B(1))$. Then the Sudakov Inequality, (\ref{3.10}), gives
     \begin{equation} 
   E\( \max_{1\le i\le n}  B(i/n)  \)\ge C\( \frac{  \log n}{n}\)^{1/2},\label{3.12}
   \end{equation}
whereas
\begin{equation}
   E\(\sup_{t\in[0,1]}B(t)\)=\sqrt{2/\pi}.
   \end{equation}
 }\end{remark}

\section{Diagonals of non-singular $M$-matrices}
 \label{sec-5}
 
\medskip	We now show that (\ref{1.34nn}) holds for a large class of non-singular $M$-matrices.  In the following we will make the assumption that    $K_{i,i}\ge  \max_{j\ne i} K_{j,i}$ and that $A=K^{-1}$ has positive row sums.  Considering the nature of the kernel of many important permanental processes this is a reasonable assumption, (see e.g.,\cite[(3.107), (3.109), and Theorem 13.1.2]{book}.

\begin{lemma}\label{lem-2.2} Let $A$ be an $n\times n$   non-singular 
$M$-matrix   with positive row sums and set $K=A^{-1}$. Assume that   $K_{i,i}\ge  \max_{j\ne i} K_{j,i}$. Then 
 \begin{equation}
  A_{i,i}\le  \frac{1}{K_{i,i}-\max_{j\ne i} K_{j,i}}.\label{6.1}
   \end{equation}
 
 \end{lemma}

\Proof   
  Using the facts that $A$ is an $M$-matrix and $   \sum_{j=1}^{n} A_{i,j}K_{j,i}=1$ and $\sum_{j\ne i}|A_{i,j}|\le A_{ii}$, we see that  
 \bea
 K_{i,i}A_{i,i}&=&1+\sum_{j\ne i}|A_{i,j}|K_{j,i}\label{2.12}\\
 &\le&1+\max_{j\ne i} K_{j,i}\sum_{j\ne i}|A_{i,j}|\le 1+\max_{j\ne i} K_{j,i}A_{i,i} \nn,
   \eea 
   which gives (\ref{6.1}).   \qed
   
   Set   \begin{equation}
  \si^{2}_{i,j}=K_{i,i} +K_{j,j}-(K_{i,j}+K_{j,i}) \quad\mbox{and}\quad    (\si_{n}^{*})^{2}=\inf_{i,j:i\ne j}   \si^{2}_{i,j}.\label{prob10}
   \end{equation}
 The fact that we can write these as squares follows from our assumption that  $K_{i,i}\ge  \max_{j\ne i} K_{j,i}$.
    
\begin{lemma}\label{lem-5.2n} Under the hypotheses of Lemma \ref{lem-2.2} assume also that   $K$ is constant along the diagonal and that 
\begin{equation}
   |K_{i,j}-K_{j,i}|\le C  \si^{2}_{i,j}\qquad\mbox{for   }C<1.
   \end{equation}
 Then
  \begin{equation}
  A_{i,i}\le  \frac{2}{ (1-C)(\si_{n}^{*})^{2}}.\label{6.1ww}
   \end{equation}

 \end{lemma}

\Proof  Consider (\ref{6.1}) and set $K_{j^{*},i}=\max_{j\ne i} K_{j,i}$. We have
\bea 
   K_{i,i}-K_{j^{*},i}&=&\frac12\( \( K_{i,i}+K_{j^{*},j^{*}}-(K_{j^{*},i} +K_{i,j^{*}}\)\)-\frac12(K_{j^{*},i}-K_{i,j^{*}})\nn\\
   &\ge &\frac12\(    \si^{2}_{i,j^{*}} -|K_{j^{*},i}-K_{i,j^{*}}|\)\\
      &\ge &\frac{(1-C)    \si^{2}_{i,j^{*}}}{2}\ge \frac{(1-C)    (\si_{n}^{*})^{2}}{2}  \nn.
 \eea
 Using this in (\ref{6.1}) we get (\ref{6.1ww}).\qed
 
 \noindent{\bf  Proof  Theorem \ref{theo-5.1qq}}  This follows from Lemma \ref{lem-5.2n} with $t_{j}=j\de/n$, $j=1,\ldots n$, for some $\de>0$ and Corollary \ref{cor-1.2} with $d_{s,t}$ replaced by $\si_{s,t}$. \qed

 \begin{remark} {\rm If $K+K^{T}$ is positive definite it is easy to see that $ \si _{i,j}$ is a metric on $\{1,\ldots,n\}$, because we can define an $n$-dimensional mean zero Gaussian random variable $\{X_{i},i\in \{1,\ldots,n\}\}$ with covariance $(K+K^{T})/2$ and 
\begin{equation}
   \si _{i,j}=(E(X_{i}-X_{j})^{2})^{1/2} \quad\mbox{and}\quad     \si_{n}^{*} =\inf_{i,j:i\ne j}  (E(X_{i}-X_{j})^{2})^{1/2}.
\end{equation}

 }\end{remark}

 We can remove the assumption that the kernel is constant on the diagonal.

  \begin{lemma}\label{lem-5.1ee} Let $A$ be an $n\times n$   non-singular 
$M$-matrix   with positive row sums and set $K=A^{-1}$. Assume that    $K_{i,i}>  \max_{j\ne i} K_{j,i}$.    Choose 
  $r_{i}=K_{i,i}/\wh K$ for some constant $\wh K$. Set  
\begin{equation}
  \wh \si_{i,j}^{2} =2\wh K- \frac{K_{i,j}}{r_{j}}-  \frac{K_{j,i}}{r_{i}}
   \end{equation}
   and assume that 
        \begin{equation}
\Bigg |   \frac{K_{i,j}}{r_{j}}-  \frac{K_{j,i}}{r_{i}} \Bigg|\le C\,\wh \si_{i,j}^{2} , \qquad C<1,\label{4.20ww}
 \ee
 for all $i,j$. Set
 \begin{equation}
   (\wh \si_{n}^{*})^{2}=\inf_{i,j:i\ne j} \wh  \si^{2}_{i ,j }.
   \end{equation}
 Then
  \begin{equation}
 r_{i} A_{i,i}\le  \frac{2}{ (1-C)(\wh \si_{n}^{*})^{2}}.\label{6.1wwj}
   \end{equation}
 \end{lemma}
 
 \Proof   Let $X=(X_{1}, \ldots, X_{n})$ be the $\al$-permanental  vector with kernel $K$.  Then $Y=(Y_{1}/r_{1}, \ldots, Y_{n}/r_{n})$  is the $\al$-permanental  vector with kernel $K_{Y}=:KR^{-1}$, where $R=\mbox{diag} (r_{1}, \ldots, r_{n})$.
 It follows from the assumption that $K_{i,i}>  \max_{j\ne i} K_{j,i}$ that this also holds for $K_{Y}$.
 Let $A_{Y}=K^{-1}_{Y}=RA$ . We see  that $A_{Y}$ is a non-singular 
$M$-matrix   with positive row sums. Consequently, (\ref{6.1wwj})  follows from Lemma \ref{lem-5.2n}.\qed

  We have the following generalization of Theorem \ref{theo-5.1qq}:
 
 \begin{theorem} \label{theo-5.1ee}
  Let $u$ be the potential density of a  transient Markov process
  in $R^{1}$.  Set 
\begin{equation}
   \wh\si_{s,t}^{2} =2 - \frac{u(s,t)}{u(t,t)}- \frac{u(t,s)}{u(s,s)}
   \end{equation}
   and assume that 
        \begin{equation}
\Big | \frac{u(s,t)}{u(t,t)}- \frac{u(t,s)}{u(s,s)}\Big|\le C\wh\si_{s,t}^{2} , \qquad C<1,\label{4.20wwef}
 \ee
 for all $|s|,|t|$ sufficiently small.    Then
    \begin{equation}
\limsup_{n\to \ff} \(\sup_{(t_{1},\ldots,t_{n})}  (\wh\si_{n}^{ *})^{2}(t_{1},\ldots,t_{n})\)\log n=\ff, 
\label{sudxx}
 \end{equation}
  implies that 
  \begin{equation}
   \sup_{t}\frac{Y_{t}}{u(t,t)}=\ff\qquad a.s.
   \end{equation} where $Y_{t}$ is the $\al$-permanental  process with kernel  $u$.
 \end{theorem}
 
 \Proof The proof is the same as the proof of Theorem \ref{theo-5.1qq}.\qed
 
 In Lemma \ref{lem-2.2} we assume that    $K_{i,i}>  \max_{j\ne i} K_{j,i}$. The following example shows that we can still get an inequality like (\ref{6.1ww}) when this condition does not hold.

   	\begin{example}\label{ex-2.1} {\rm Consider the covariance matrix $\BB$ of $(B(1),\ldots,B(n))$, where $\{B(t),t\in R_ { +}\}$ is standard Brownian motion. Obviously, $\BB_{i,i}-\max_{j\ne i} \BB_{i,j}=0$. However, $\BB^{-1}$ is a tri-diagonal matrix with all diagonal elements equal  to 2, except that $(\BB^{-1})_{n,n}=1$ and all  off diagonal elements that are not zero equal to -1. In this case
\begin{equation}
   (\BB^{-1})_{i,i}\le 2=\frac{2}{(\si_{n}^{*})^{2}}.
   \end{equation}
   (Here $(\si_{n}^{*})^{2}=\min_{i\ne j}E(B(i)-B(j))^{2}=1$.)

We can use this to create another interesting example. Let $D$ be a diagonal matrix with entries $1,\ldots,n$. Let $\wt\BB=D^{-1}\BB  $. This matrix has entries 1 on and above the diagonal and $ \wt\BB_{ i,j}=j/i$ for $1<j<i$. The diagonal entries of $\wt\AA=(\wt\BB)^{-1}$ are $\wt\AA_{i,i}=2 B_{i,i}=2i$, $1\le i\le n-1$ and $\wt\AA_{n,n}= B_{n,n}=n$.  Set
\begin{equation}
   \phi^{2}_{i,j}= \wt\BB_{ i,i}+ \wt\BB_{ j,j}- \wt\BB_{ i,j} -\wt\BB_{ j,i}
   \end{equation}
   and 
   \begin{equation}
  ( \phi_{n}^{*})^{2} =\min _{\stackrel{1\le i,j\le n}{i\ne j}}   \phi^{2}_{i,j}. 
   \end{equation}
   The minimum is achieved at $\phi^{2}_{n ,n-1}=1/n$. Therefore we have
   \begin{equation}
   (\AA)_{i,i} \le \frac{2}{(\phi_{n}^{*})^{2}}=2n\qquad 1\le i\le n.\label{2.54}
   \end{equation}
   The maximum on the left-hand side of (\ref{2.54}) is $ (\AA)_{n-1,n-1}=2(n-1)$, since $ (\AA)_{n,n}=n .$ 
   
 }\end{example}

\begin{lemma} \label{lem-5.4}
When  $u$ is the potential density of a    transient Markov process
  in $R^{1}$, (\ref{4.20wwef}) always holds with $C=1$.
  \el
  
  \Proof      We need to show that 
   \begin{equation}
\Big | \frac{u(s,t)}{u(t,t)}- \frac{u(t,s)}{u(s,s)}\Big|\le  2 - \frac{u(s,t)}{u(t,t)}- \frac{u(t,s)}{u(s,s)},\label{5.20}
 \ee 
  Without loss of generality we assume that  $u(s,t)/u(t,t)\ge u(t,s)/u(s,s)$. Then (\ref{5.20}) is equivalent to
  \begin{equation}
   \frac{u(s,t)}{u(t,t)}\le 1.
   \end{equation}
 It follows from \cite[Lemma 3.4.3]{book} that when  $u$ is the potential density of a    transient Markov process, in $R^{1}$,     this always holds.\qed

    \section{Permanental processes with a kernel that is the potential density of  a L\'evy process}\label{sec-6}

  \noindent{\bf Proof of  Theorem  \ref{theo-new} }  It follows from Lemma \ref{lem-5.1} that (\ref{3.25nn}) is the same as (\ref{4.20wwe}). Therefore (\ref{5.7qq}) follows from (\ref{sudx}) with $t_{1},\ldots,t_{n}$ replaced by $\de/n,2\de/n,\ldots,\de$. \qed   
  
  The next lemma is used in the proof of Theorem \ref{theo-1.6}
    
\begin{lemma} \label{lem-6.1}
Suppose that  $\ell$  and h are positive,   quasi-monotonic slowly varying functions (see \cite[Section 2.7]{BGT}) at infinity. Set 
 \bea 
  | \II  (   \la)| =  \frac{\ell(|\la|)}{|\la|}\quad\mbox{and}\quad   | \RR  (   \la)| =  \frac{h(|\la|)}{|\la|}.\label{newdef276}
 \eea
 If $\RR\in L^{1}$ and
\begin{equation}
\int_{1/z}^{\ff} \RR( \la) \,d \la  \ge B   \ell( 1/|z|),\label{8.16x}
\end{equation}
as $|z|\to 0$  with $   B>\frac{\pi}{2}$, then
   \begin{equation}
 \Big|  \int_{0}^\ff   \sin(\la z)\,  \II  (   \la) \,d\la\Big|\le   C\int_{0}^\ff (1- \cos(\la z) )\RR ( \la) \,d\la\label{8.0}
   \end{equation}
   for some $C<1$, and all $|z|$ sufficiently small. 
   Furthermore,  
  \be   \int_{0}^\ff (1-\cos (\la z))\RR  (\la)\,d\la\sim  \int_{1/|z|}^{ \ff}\RR ( \la) \,d\la\label{6.5}
  \ee as $|z|\to 0$.
 \end{lemma}

 \Proof  It suffices to show (\ref{8.0}) for $z>0$. 
       By \cite[(1.43)]{S} 
\begin{equation}
\int_{0}^\ff  \frac{ 1_{\{\la z\leq 1\}}-e^{i\la z} }{\la}\ell(   \la)\,d\la\sim \ell( 1/z)\int_{0}^\ff    \frac{ 1_{\{\la z\leq 1\}}-e^{i\la z} }{\la}\,d\la,\label{8.10}
\end{equation}
as  $z\to0$. 
Taking the imaginary part of (\ref{8.10}) we see that 
\begin{equation}
 \int_{0}^\ff   \sin(\la z)\,  \II  (   \la) \,d\la= \int_{0}^\ff  \,{   \sin(\la z)\over \la}\ell(   \la)\,d\la\sim  \ell( 1/z)\int_{0}^\ff  \,{   \sin(s)\over s}\,ds ,\label{8.11}
\end{equation}
 as $z\to 0$.
  Therefore, by \cite[3.721]{GR},   
  \begin{equation}
 \int_{0}^\ff   \sin(\la z)\,  \II  (   \la) \,d\la \sim \frac{\pi  }{2 }  \ell( 1/z),\label{8.6}
  \end{equation}
  as $z\to 0$.
  
To use below we note that by a change of variables
 \bea
 \int_{0}^\ff  \,{ 1_{\{\la z\leq 1\}}- \cos(\la z)\over \la}\,d\la&=&\int_{0}^\ff  \,{ 1_{\{s\leq 1\}}- \cos(s)\over s}\,ds\label{6.9}\\
 &=&2\int_{0}^1{ \sin^{2}(s/2)\over s}\,ds- \int_1^\ff  {\cos(s)\over s}\,ds.\nn
 \eea
  Therefore the first integral in (\ref{6.9}) is a constant which we denote by $c_{0}$. It is easy to see that  $c_{0}<\ff$. 
The first of the  last two integrals in (\ref{6.9}) is bounded by 1/4, and by integration by parts, that the second of these last two integrals is bounded by 2.  
 
By (\ref{8.10}), (\ref{6.9})  and \cite[Theorem 14.7.2]{book}
\begin{eqnarray}
&&\int_{0}^\ff (1- \cos(\la z) )\RR( \la) \,d\la
\label{8.15}\\
&&\qquad=\int_{0}^{\ff} (1_{\{\la z\leq 1\}}- \cos(\la z) ) \RR( \la) \,d\la+\int_{1/z}^{\ff}  \RR( \la) \,d\la   \nonumber\\
&& \qquad\sim   {\it h(1/z)} \int_{0}^\ff  \,{ 1_{\{\la z\leq 1\}}- \cos(\la z)\over \la}\,d\la+\int_{1/z}^{\ff}  \RR( \la) \,d\la\nn\\
&&\qquad = c_{o} {\it h(1/z)}+\int_{1/z}^{\ff}  \RR( \la) \,d\la\sim \int_{1/z}^{\ff}  \RR( \la) \,d\la\nn.
\end{eqnarray}
as $z\to 0$.    Thus we obtain (\ref{8.0})   and also (\ref{6.5}).  (See (\ref{1.28qq}).)\qed

  \medskip	
  \noindent{\bf Proof of Theorem \ref{theo-1.6} }  The characteristic exponent of this process 
        \bea
   \psi(\la)&=&-\int_{\ff}^{\ff} \(e^{i\la x}-1- i\la x 1_{\{|x|<1\}}\)\,    \nu(dx)  \label{6.11q}\\
   &\sim&\frac{\pi}{2} |\la| g(\la)+i(p-q) \la \int_{1/ \la}^{1} \frac{g(1/x)}{x}\,dx \nn\\
      &\sim&\frac{\pi}{2}|\la| g(\la)+i(p-q) \la \int_{  1}^{\la} \frac{g(s)}{s}\,ds \nn 
    \eea
 as $\la\to\ff$.   Note that  the   $Re\,\psi(\la)=o(Im\,\psi(\la))$ as $\la\to\ff$.
   We show how to obtain (\ref{6.11q}) in Section \ref{subsec-7.2}.

  We first consider the case when $p\ne q$. It follows from (\ref{6.11q}) that   
    \be 
   \II_{\bb}(\la)\sim \frac{1}{|p-q| \la\int_{  1}^{\la} \frac{g(s)}{s}\,ds}:=\frac{\ell(\la)}{\la}\label{6.12}
 \ee
 and  
\be 
   \RR_{\bb}(\la)\sim \frac{(\pi/2) g(\la)}{|p-q|^{2}  \la \(\int_{  1}^{\la} \frac{g(s)}{s}\,ds\)^{2}}.\label{6.13}
 \ee
  as $\la\to\ff$.
Note that  
\begin{equation}
\frac{(\pi/2) g(\la)}{|p-q|^{2}  \la \(\int_{  1}^{\la} \frac{g(s)}{s}\,ds\)^{2}}=-\frac{(\pi/2)}{|p-q|}  \frac{d}{d\la}\ell(\la)
   \end{equation}
 which implies that
 \begin{equation}
  \int_{1/z}^{\ff}     \RR_{\bb}(\la)\, d\la\sim \frac{(\pi/2)}{|p-q|}\ell(1/z).\label{6.16}
   \end{equation}
  Comparing this with (\ref{6.12})
we see that  (\ref{8.16x}) holds for all $B< (\pi/2)/|p-q|$. Obviously, we can take $B>(\pi/2)$ as long as $p\ne q$. Also,  
 by (\ref{1.28qq}) and  (\ref{6.5}),   
 \begin{equation}
   \si^{2}(z)\sim \frac{1}{|p-q|^{2}}\(  \int_{1}^{1/z}    \frac{g(s)}s\,ds \)^{-1}.
   \end{equation}
Therefore, (\ref{6.19s}) follows from Theorem \ref{theo-new}. 

Note that we require that $\RR_{\bb}\in L^{1}(R_{+})$. That is why we impose the condition in (\ref{1.30}).

\medskip	When $p=q$, $\psi(\la)$ is real and symmetric and 
\begin{equation}
   \RR_{\bb}(\la)\sim \frac{2}{\pi |\la| g(\la)}.
   \end{equation}
Conditon (\ref{8.16x}) in Lemma \ref{lem-6.1} is trivially satisfied and, by 
 by (\ref{1.28qq}) and  (\ref{6.5})     
  \be \si^{2}(z)\sim \frac{4}{\pi^{2}}\int_{1/|z|}^{ \ff}\frac{1}{ |\la| g(\la)} \,d\la\label{6.5m}
  \ee as $|z|\to 0$.
Therefore, (\ref{1.33}) follows from Theorem \ref{theo-new}. 
\qed

 \noindent{\bf Details for (\ref{1.31}) } 
 This is simple for symmetric processes, so we only need to check (\ref{1.31}) for $p\ne q$. By    (\ref{1.28qq}), for all L\'evy processes,  
 \begin{equation}
   R_{\bb}(z)=u^{\bb}(0)-\frac{\si^{2}(z)}{2}.\label{6.17}
   \end{equation}
 For the processes with L\'evy measure given by (\ref{1.30s}) we see by (\ref{8.6})   that as $z\to 0$
 \begin{equation}
   H_{\bb}(z)\sim \frac{\ell(1/z)}{2}.
   \end{equation}
 In addition by (\ref{1.28qq}), (\ref{8.15}) and (\ref{6.16})   
 \begin{equation}
   \si^{2}(z)\sim \frac{2}{\pi}\int_{1/z}^{\ff}R(\la)\,d\la\sim \frac{1}{|p-q|}\ell(1/z).
   \end{equation}
 Therefore,  
 \begin{equation}
      H_{\bb}(z)\sim \frac{|p-q|}{2}  \si^{2}(z)\quad\mbox{and}\quad   H_{\bb}(-z)\sim -\frac{|p-q|}{2}  \si^{2}(z) \label{6.20}.
   \end{equation}
 Adding (\ref{6.17}) and (\ref{6.20}) we get (\ref{1.31}).

 \medskip	 \noindent{\bf Proof of  Corollary  \ref{cor-1.4} } To show that (\ref{3.25nn}) is satisfied it suffices to show that  for all $z>0$ sufficiently small 
   \begin{equation}
 \Big|\int_{0}^\ff   \sin(\la z)\,  \II  _{\bb}(   \la) \,d\la \Big| \le   C\int_{0}^\ff (1- \cos(\la z) )\RR_{ \bb}( \la) \,d\la\label{3.26m}
   \end{equation}
   for some $C<1$.
   To simplify the proof we  assume that $ \II  _{\bb}(   \la)\geq 0 $ and take $ \II  _{\bb}(   \la) $ and $ \RR  _{\bb}(   \la) $ to be  non-increasing functions. It is easy to extend the proof  to the case in which $| \II  _{\bb}(   \la) |$ and $ \RR  _{\bb}(   \la) $ are asymptotic to  non-increasing functions as $\la\to\ff$. We have  
     \begin{equation}
     \int_{0}^\ff   \sin(\la z)\,   \II  _{\bb}(   \la) \,d\la  \le    \int_{0}^ {\pi/z}  \la z\,   \II  _{\bb}(   \la)\,d\la\label{5.12}
   \end{equation}  
   because, since $\II_{\bb}(\la)$ is decreasing, for all $k\ge 1$
   \begin{equation}
      - \int_{(2k-1)\pi/z}^{(2k )\pi/z}  \sin(\la z)\,  \II  _{\bb}(   \la) \,d\la\ge     
      \int_{(2k)\pi/z}^{(2k+1 )\pi/z}  \sin(\la z)\,  \II  _{\bb}(   \la) \,d\la.
   \end{equation}
 Also
     \bea 
     \frac{1}{2} \int_{0}^\ff (1- \cos(\la z) )\RR_{ \bb}( \la) \,d\la&=&
 \int_{0}^{\ff}\sin^{2}\(\frac{\la z}{2}\) \RR_{\bb}(\la)\,d\la \nn\\&\ge&   \int_{\pi/(2z)}^{\ff}\sin^{2}\(\frac{\la z}{2}\) \RR_{\bb}  (\la)\,d\la\label{}\\
 &=& \sum_{k=0}^{ \ff}  \int_{\pi(1+4k)/(2z)}^{\pi(1+4(k+1))/(2z)}\sin^{2}\(\frac{\la z}{2}\) \RR_{\bb}  (\la)\,d\la\nn\\
 &\ge& \sum_{k=0}^{ \ff}  \int_{\pi(1+4k)/(2z)}^{\pi(3+4k)/(2z)}\sin^{2}\(\frac{\la z}{2}\) \RR_{\bb}  (\la)\,d\la\nn.
 \eea
 Note that if $\pi(1+4k)/(2z)\leq\la\leq \pi(3+4k)/(2z)$ then
\be \pi/4+k\pi\le \la z/2\le 3\pi/4+k\pi,\ee
 and consequently $\sin^{2}(\la z/2)\ge 1/2$.  Therefore, 
 \begin{equation}
  \int_{\pi(1+4k)/(2z)}^{\pi(3+4k)/(2z)}\sin^{2}\(\frac{\la z}{2}\) \RR_{\bb}  (\la)\,d\la\ge \frac{1}{2} \int_{\pi(1+4k)/(2z)}^{\pi(3+4k)/(2z)} \RR_{\bb}  (\la)\,d\la. \label{wwww}
 \end{equation}
  Furthermore, since $\RR_{\bb}(\la)$ is decreasing,  for all   $z> 0$ sufficiently small  
 \bea
 &&   \frac{1}{2} \int_{\pi(1+4k)/(2z)}^{\pi(3+4k)/(2z)} \RR_{\bb}  (\la)\,d\la\\
&&\qquad\ge  \frac{1 }{4}  \int_{\pi(1+4k)/(2z)}^{\pi(3+4k)/(2z)} \RR_{\bb}  (\la)\,d\la+\frac{1 }{4}  \int_{\pi(3+4k)/(2z)}^{\pi(5+4k)/(2z)} \RR_{\bb}  (\la)\,d\la \nn\\
    &&\qquad= \frac{1 }{4}    \int_{\pi(1+4k)/(2z)}^{\pi(1+4(k+1))/(2z)} \RR_{\bb}  (\la)\,d\la.\nn
   \eea 
   Putting all this together we see that  for all $z $ sufficiently small   
      \be 
 \int_{0}^\ff (1- \cos(\la z) )\RR_{ \bb}( \la) \,d\la  
 \ge  \frac{1 }{2}    \int_{\pi/(2z)}^{\ff}  \RR_{\bb}  (\la)\,d\la\label{3.41q} .
   \ee 
   Combining (\ref{5.12}) and (\ref{3.41q}) and using the hypothesis  (\ref{5.11}) we get (\ref{3.26m}) for   some $C<1$. 
   
   We see from (\ref{1.28-}) and (\ref{3.41q}) that
   \begin{equation}
   \si^{2}(z)\ge C \int_{\pi/(2z)}^{\ff}  \RR_{\bb}  (\la)\,d\la
   \end{equation}
  which implies by Theorem \ref{theo-new} that  if
   \begin{equation}
 \limsup_{n\to \ff} \(  \int_{\pi n/(2 \de)}^\ff  \RR_{ \bb}( \la) \,d\la\)\log n=\ff,  \label{5.7ww}
   \end{equation}
  then the $\al$-permanental  process with kernel $u^{\bb}$ is  unbounded  almost surely.  It is easy to see that, by interpolation, this is equivalent to (\ref{5.7q}).\qed

 \begin{remark} {\rm \label{rem-5.3}   We consider (\ref{5.11}) for   $\II_{\bb}(\la)$ and $\RR_{\bb}(\la)$ assymptotic to $\II(\la)$ and $\RR (\la)$ as $\la\to\ff$; (see (\ref{6.12}) and (\ref{6.13})). In this case Corollary  \ref{cor-1.4} it is not much cruder than the estimates given in   Theorem \ref{theo-1.6}.  Since   in this case $\la\II_{\bb}(\la)$ is slowly varying at infinity we see that the left-hand side of (\ref{5.11}) is asymptotic to $\pi\ell(1/|z|)$ as $|z|\to 0$. By (\ref{8.16x}) and the fact that $\ell$ is slowly varying the  right-hand side of (\ref{5.11}) is asymptotic to $ C\pi\ell(1/|z|)/(2|p-q|)$ as $|z|\to 0$. Therefore  (\ref{5.11})  holds for $C>2|p-q|$. Therefore   Corollary  \ref{cor-1.4} gives the results  obtained in Example \ref{ex-1.1}   when $p\ne q$ and $|p-q|<1/2$.  
 }\end{remark}

  \section{Appendix}\label{sec-App}
  
  \subsection{A property of the potential density of a transient Markov process}\label{subsec-7.1}\bl
  If $\{u(s,t),s,t\in T\}$ is the potential density of a transient Markov process $X_{t}$ with state space $T$, 
  then for all $(t_{1},\ldots,t _{n})$ in $T$, the matrix $K$ with elements $\{u(t_{i},t_{j})\}_{i,j=1}^{n}$ is invertible and  its inverse  is a non-singular $M$-matrix.
  \el

\Proof    This proof is a portion of the proof of \cite[Theorem 13.1.2]{book}. We define the following stopping
time:
\begin{equation}
\si=\inf \{ t\geq 0\,|\, X_{ t}\in \{ t_{ 1},\ldots, t_{ n}\}\cap \{X_{0}
\}^{ c}\}\label{ap.1}
\end{equation} (note that $\si$ may be infinite). Let $\{ L^{ x}_{ t}\,;(
x,t)\in S\times R^{ 1}\}$ be the local times of $X$. Since
$u(t_{ i},t_{ j} )$ is the $0$-potential 
density of $X$, we  can normalize the
local time so  that
\be
u(t_{ i},t_{ j} )=E^{ t_{ i}}\( L^{ t_{ j}}_{\ff} \).\label{asp.0.49}
\ee

      Using (\ref{asp.0.49}) and  the  strong Markov property, we see that
\begin{eqnarray} u(t_{ i},t_{ j} ) &=& E^{ t_{ i}}\( L^{ t_{ j}}_{\si}
\)+E^{ t_{ i}}\(E^{ X_{\si}}( L^{ t_{ j}}_{\ff});\,\si<\ff
\)\label{ap.2}\\   &=& d_{ i,j}+
\sum_{ k=1}^{ n}h_{ i,k}u(t_{ k},t_{ j} ),
\nonumber
\end{eqnarray}
      where $d_{ i,j}=E^{ t_{ i}}( L^{ t_{ j}}_{\si} )$ and
$h_{ i,k}=P^{ t_{ i}}(X_{\si}= t_{ k})$.

Let  $D=\{ d_{ i,j}\}_{ 1\leq i,j\leq n}$ and $H=\{ h_{ i,j}\}_{ 1\leq i,j\leq
n}$.  We can write (\ref{ap.2}) as
\begin{equation}
K=D+HK\label{ap.3}
\end{equation}
      so that  $( I-H)K=D$. Moreover, $D$ is a diagonal matrix with all its
diagonal elements strictly positive. This follows  because, starting from
$X_0=t_i$,
$\si>0$, which implies that  each $b_{i,i}>0$. On  the other hand, the
process is killed the first time it hits any $t_j\ne t_i$. Thus,  starting from
$t_i$,
$L_\si^{t_j}=0$,
$j\ne i$.

Since $D$ is invertible, both  $( I-H)$ and  $ K$ are invertible and
\begin{equation}
      {K}^{ -1}=D^{ -1}( I-H).\label{ap.4}
\end{equation}
      It is clear that $H\geq 0$. It follows from this that $K^{ -1}$ has
negative off diagonal elements.  Moreover, since $h_{i,i}=0$ it follows that $K^{ -1}$ has
positive diagonal elements. Therefore, $K$ is a non-singular M-matrix. Furthermore,
\be
\sum_{ j=1}^{ n}h_{ i,j}=P^{ t_{i}}(\si<\ff)\leq 1\qquad \forall \,i=1\ldots,n,
\ee from which it follows that $K^{ -1}$  has positive row sums.
\qed

  \subsection{Derivation of (\ref{6.11q})}\label{subsec-7.2}

 We have
 \be 
 {Im}\, \psi(\la)=-(p-q)\int_{0}^{\ff} \(\sin\la x-\la x 1_{\{|x|<1\}}\)\nu(dx).\label{7.7}
 \ee 
Let  $\nu_{1}(dx):= (p-q) \nu( x)$. Then $   {Im}\, \psi(\la)$  is equal to
 \be  
   -\int_{0}^{1/\la} \(\sin\la x-\la x \)\nu_{1}(dx) 
   +\la \int_{1/\la}^{1}   x \nu_{1}(dx) -\int_{1/\la}^{\ff}  \sin\la x \nu_{1}(dx) .\label{7.8}
 \ee 
 Using $|\sin\la x-\la x|\le |\la x|^{3}$ in the first of these integral and $|\sin\la  |\le  1$ in the the last  of these integral we see that their absolute values   are both $O(g(\la)/\la)$ as $\la\to \ff$.  Consequently
 \begin{equation}
    {Im}\, \psi(\la)\sim (p-q)\la \int_{1/\la}^{1}   \frac{g(1/x)}{x} \,dx =(p-q)\la \int_{1}^{\la}   \frac{g(s)}{s} \,ds
   \end{equation} as $\la\to\ff$.
 
\medskip	The asymptotic behavior of  $Re\,\psi(\la)$ as $\la\to\ff$  follows from the next lemma.  

 	\bl\label{lem-7.2} Let $g(\cd)$ be a slowly varying function at infinity. Then
 \begin{equation}
   \int_{0}^{\ff}(1-\cos \la x)\frac{g(1/x)}{x^{2}}\,dx\sim \frac{\pi}{2}  \la  g(\la)   \label{7.9}
   \end{equation}
 as $\la\to\ff$.  
 \el
 
\Proof We write the left-hand side of (\ref{7.9}) as 
 \begin{equation}
    \la g(\la)  \int_{0}^{\ff}\frac{(1-\cos s)}{s^{2}}\frac{g(\la/s)}{g(\la)}\,ds 
    \label{7.10}.
   \end{equation}
 Consider
 \begin{equation}
   \int_{0}^{M}\frac{(1-\cos s)}{s^{2}}\frac{g(\la/s)}{g(\la)}\,ds+ \frac{1}{g(\la)}\int_{M}^{\ff}\frac{(1-\cos s)}{s^{2}}g(\la/s)\,ds\label{7.11}
   \end{equation}
   Note that by \cite[Theorem 1.5.6]{BGT}, for $s\in (0,M]$, $g(\la/s)/g(\la)\le C(s ^{ -\ep}\vee 1)$, for any $\ep>0$, and some constant depending on $M$ and $\ep$. Therefore, by the dominated convergence theorem 
 we see that the limit, as $\la\to\ff$ of the first integral in (\ref{7.11}) is 
\be
   \int_{0}^{M}\frac{(1-\cos s)}{s^{2}} \,ds.\label{7.13}
   \ee

The second integral in (\ref{7.11}) is bounded by  
 \begin{equation}
     \frac{2}{\la g(\la)}  \int_{M/\la}^{\ff}\frac{g(1/s) }{s^{2}}\,ds=     \frac{2}{\la g(\la)}  \int_{0}^{\la/M} g(v)\,dv\label{7.14q}.
   \end{equation}
   We need a condition on $g$ near 0. This is given implicitly by the statement that $\nu$ is a L\'evy measure, which requires that
    \begin{equation}
   \int_{1}^{\ff}\frac{g(1/|x|)}{x^{2}}\,dx=    \int _{0}^{1} g(v)  \,dv \le C'<\ff
   \end{equation}
  Therefore, since $g$ is slowly varying at infinity, the second integral in (\ref{7.14q}) is bounded by $3 g(\la/M) /(Mg(\la))$
 which goes to $3/M$ as $\la\to\ff$. Therefore,  taking $\la\to\ff$, we see that (\ref{7.11}) is equal to 
 \begin{equation}
   \int_{0}^{M}\frac{(1-\cos s)}{s^{2}} \,ds+O(1/M)     \end{equation}
  for all $M$. This gives us (\ref{7.9}) because, by integration by parts,    
  \begin{equation}
    \int_{0}^{\ff}\frac{(1-\cos s)}{s^{2}} \,ds =    \int_{0}^{\ff}\frac{ \sin s }{s } \,ds={\pi \over 2}, 
   \end{equation}
 by \cite[3.721]{GR}.
 \qed

  \vspace{1 in} 
\def\noopsort#1{} \def\printfirst#1#2{#1}
\def\singleletter#1{#1}
            \def\switchargs#1#2{#2#1}

\def\bibsameauth{\leavevmode\vrule height .1ex
            depth 0pt width 2.3em\relax\,}
\makeatletter
\renewcommand{\@biblabel}[1]{\hfill#1.}\makeatother
\newcommand{\bysame}{\leavevmode\hbox to3em{\hrulefill}\,}

   \end{document}